\documentclass{article}
\usepackage{epsf}
\usepackage{latexsym}
\usepackage{amsfonts}
\usepackage{longtable}

\begin{document}

\title{An Alternative Computational Approach to the $3n+1$ Conjecture}

\author{Kevin P. Thompson}

\date{}

\thispagestyle{empty}
\renewcommand\thispagestyle[1]{} % Prevents \maketitle from inserting page no.

\maketitle
\begin{abstract}
An alternative computational approach to the $3n+1$ conjecture is presented that may be theoretically capable of confirming the conjecture.
\end{abstract}

The $3n+1$ conjecture, also known as the Collatz Conjecture, the Ulam conjecture, or the Syracuse problem, was first proposed in 1928 by Lothar Collatz \cite{Lagarias}. The conjecture states that the following sequence of numbers always reaches 1: beginning with a natural number $n$, if $n$ is even, divide it by 2; if $n$ is odd, multiply it by 3 and add 1. Repeat the process with the resulting number.

The primary computational approach to the conjecture has been to attempt to verify that the sequence reaches 1 for the first $n$ natural numbers. At the time of writing, the record for this verification stands in excess of $2.7\times 10^{16}$ \cite{Oliveira}. However, if the conjecture is true, then this computational approach is in vain. It is only capable of disproving the conjecture or providing some confidence that it might be true.

Since the current approach has thus far failed to disprove the conjecture, it would be helpful to attempt a computational approach that, at least in theory, may have the potential to prove the conjecture. This paper sets forth one such strategy.

\section{Motivation}

Consider a proof of the $3n+1$ conjecture by mathematical induction. The sequence obviously reaches 1 for the first several natural numbers. Therefore, let $n$ be a natural number and assume the sequence for $n-1$ reaches 1. Our objective is to show that the sequence for $n$ also reaches 1. But, if any element of the sequence for $n$ appears in any previous sequence, then the sequence for $n$ must also reach 1 by assumption. Therefore, at a minimum, we must show that the sequence for $n$ contains an element in a previous sequence.

At this point the proof breaks down into many, many cases. But, we can take a systematic approach to finding subsets of real numbers that can be described as having a sequence that reaches 1. We will consider congruence classes mod $a$ where $a$ is even since members of these classes can be described as odd or even.

If we examine congruence classes mod 2, we see that the second element of the sequence for numbers in $\overline{0}_2$ appears in a previous sequence. If $n=2k$ for some natural number $k$, then $n$ is even and the next element of the sequence is $k$. By assumption, the sequence for $k<n$ reaches 1 so the sequence for $n$ must also reach 1. The sequence for $\overline{1}_2$ follows the pattern $(2k+1,6k+4,3k+2,\dots)$ at which point the 3rd element may be either even of odd. Therefore, no general statement can be made concerning $\overline{1}_2$. For these congruence classes, the sequence for 50\% of the classes reaches 1.

Now consider congruence classes mod 4. Even numbers are accounted for since the congruence classes $\overline{0}_4$ and $\overline{2}_4$ together are equivalent to $\overline{0}_2$. For members of $\overline{3}_4$, the elements of the sequence follow the pattern
\begin{center}
$4k-3$ (odd)

$12k-8$ (even)

$6k-4$ (even)

$3k-2$
\end{center}
Since the fourth element of the sequence is less than the sequence anchor, the sequence reaches 1 by assumption. The sequence pattern for congruence class $\overline{1}_4$ is inconclusive for the same reason as $\overline{1}_2$ seen above. For these congruence classes, the sequence pattern for 75\% of the congruence classes reaches 1.

Progressing to congruence classes mod 6, we again skip $\overline{a}_6$ when $a$ is even since these are even numbers accounted for by $\overline{0}_2$. The only successful congruence class is $\overline{1}_6$. This can be seen by considering the sequence pattern for $\overline{1}_4$:
\begin{center}
$4k-1$ (odd)

$12k-2$ (even)

$6k-1$
\end{center}
The sequence anchor of sequences for members of $\overline{1}_6$ appears as the 3rd element of the previous sequence pattern for members of $\overline{1}_4$. Therefore, sequences for members of $\overline{1}_6$ reach 1 by assumption.

To compute our coverage (or percent complete), we must first realize that $\overline{1}_6$ and $\overline{3}_4$ overlap. We must examine the corresponding congruence classes mod 12, the least common multiple (LCM) of the moduli 2, 4, and 6. Mapping the successful congruence classes, $\overline{0}_2$ is equivalent to $\overline{0}_{12}$, $\overline{2}_{12}$, $\overline{4}_{12}$, $\overline{6}_{12}$, $\overline{8}_{12}$, and $\overline{10}_{12}$; $\overline{3}_4$ is equivalent to $\overline{3}_{12}$, $\overline{7}_{12}$, and $\overline{11}_{12}$; and, $\overline{1}_6$ is equivalent to $\overline{1}_{12}$ and $\overline{7}_{12}$. The only congruence classes mod 12 not accounted for are $\overline{5}_{12}$ and $\overline{9}_{12}$. Therefore, the sequence pattern for approximately 83\% of the congruence classes reaches 1. 

The systematic approach can be continued and implemented as an algorithm. The cases need not be so simple, but there are only a few successful scenario types. As another example, the 12th element of sequences for members of $\overline{33}_{128}$ is the 9th element of sequences for members of $\overline{25}_{96}$:

\begin{center}
\begin{tabular}{ c c }
 $128k-33$ & $96k-25$ \\
 $384k-98$ & $288k-74$ \\
 $192k-49$ & $144k-37$ \\
 $576k-146$ & $432k-110$ \\
 $288k-73$ & $216k-55$ \\
 $864k-218$ & $648k-164$ \\
 $432k-109$ & $324k-82$ \\
 $1296k-326$ & $162k-41$ \\
 $648k-163$ & $486k-122$ \\
 $1944k-488$ &  \\
 $972k-244$ &  \\
 $486k-122$ &  \\
\end{tabular}
\end{center}

Cases such as the ones illustrated above abound. Each case illustrates a congruence class mod $a$ of integers whose members have a $3n+1$ sequence that reaches 1 based on the mathematical induction assumption. In order for the proof to succeed, enough successful congruence classes (of various moduli) will need to be found so that the congruence classes modulo the LCM of the individual moduli all reach 1.

\section{The Computational Approach}

With the even numbers accounted for by $\overline{0}_2$, and with $\overline{3}_4$ being the only successful congruence class mod 4, we wish to identify as many congruence classes $\overline{c}_b$ with $b$ even, $b\ge 6$, $c$ odd, and $1\le c<b$, whose $3n+1$ sequence can be proven to contain a member of a previous sequence. In terms of sequence patterns, this is equivalent to identifying patterns $n=bk-c$ whose sequence can be proven to contain a member of a sequence less than $n$. This is done by looking for the sequence to directly become less than its anchor or by the sequence otherwise containing a member of a previous sequence.

Begin with $b=6$ and examine the sequence $6k-1$. Construct the sequence as far as is possible (i.e. until the factor becomes odd). Check if the sequence becomes less than the sequence anchor. If not, construct all previous sequences (i.e. $4k-1$ and $4k-3$) and look for elements of $6k-1$ in these sequences. If a sequence pattern is found to reach 1, the successful congruence classes are restated with the new modulus being the LCM of all the individual moduli to determine if a partition the integers has been found. A measure of completeness is also calculated as the percent of congruence classes found to be successful.

Repeat the process in order of increasing even numbers $b$. For each selection of $b$ and $c$, the sequence for $bk-c$ is constructed as far as is possible with checks to see if a sequence element becomes less than the anchor. Then, all previous sequences are constructed searching for overlapping elements with the given sequence.

\section{Dependencies}

The success of the algorithm is first dependent on an upper bound existing for how long it takes a $3n+1$ sequence to join a previous sequence. This condition is equivalent to an upper bound existing for the stopping times of all sequences. If such an upper bound does not exist, this algorithm will continually approach 100\% completion but never achieve it. But, if an upper bound exists, the algorithm may eventually come to completion.

This algorithm is capable of detecting at least some types of non-trivial cycles. The algorithm would enter an infinite loop if certain sequences contained a cycle.

\section{Analysis}

It is intuitively obvious that we must have $b\equiv 0\bmod (2^t)$ for some $t$ for otherwise we would not be able to perform division by 2 operations. Computationally, the only successful moduli less than 1,500 were of the form $2^t3^s$, $t\ge 1$ and $s\ge 0$. The search has been narrowed and is continuing for $bk-c$ such that $b=2^t3^s$ and $c$ is odd. 

In addition, if a sequence was found to join a previous sequence, the modulus of the previous sequence was always of the form $2^t3^s$ for the initially discovered sequences. Therefore, when looking for a previous sequence for the current sequence to join, only previous sequences with a modulus of the form $2^t3^s$ are now being considered.

If a sequence $bk-c$ joins a previous sequence $dk-e$ with $b=2d$, then it appears $c$ and $e$ must satisfy $c-e=2^t3^s$ for some $s\ge 0$ and $t\ge 1$.

There appear to be dead zones between moduli in which no successful congruence classes are found. Between two moduli $2^t$ and $2^{t+1}$, success only occurs if the modulus is less than or equal to the midpoint, $3\times 2^{t-1}$.

One of the chief concerns is of course whether the algorithm will complete. Analysis of the progress of the algorithm to date may give some indications. Percent complete gains are most significant for moduli of the form $2^t$. Therefore, looking at the change in the completion percentage after the processing of each modulus of a power of 2 may indicate whether progress will continue at a sufficient rate to reach completion.

The table below summarizes the approximate increases in the completion percentage for such moduli. Interestingly, although the general pattern is a decrease in the rate, the rates are not strictly decreasing. Sudden increases and lack of a consistent halving of the progress are positive indications the algorithm may reach completion.

\begin{center}
\begin{tabular}{ c c }
 Factor & Pct complete change \\
 $2^4$ & 4.16667\% \\
 $2^5$ & 3.47222\% \\
 $2^6$ & 0.86806\% \\
 $2^7$ & 1.30208\% \\
 $2^8$ & 1.27315\% \\
 $2^9$ & 0.53048\% \\
 $2^{10}$ & 0.42438\% \\
 $2^{11}$ & 0.31266\% \\
 $2^{12}$ & 0.20703\% \\
 $2^{13}$ & 0.29246\% \\
 $2^{14}$ & 0.15194\% \\
 $2^{15}$ & 0.11030\% \\
 $2^{16}$ & 0.18778\% \\
 $2^{17}$ & 0.08812\% \\
 $2^{18}$ & 0.07699\% \\
 $2^{19}$ & 0.07083\% \\
 $2^{20}$ & 0.05061\% \\
 $2^{21}$ & 0.08550\% \\
 $2^{22}$ & 0.04891\% \\
\end{tabular}
\end{center}

It is also interesting to examine the approximate increases in percent complete between moduli of the form $2^t$. This will give an indication of the contribution of mixed factors of the form $2^t3^s$.

\begin{center}
\begin{tabular}{ c c }
 Factor range & Pct complete change \\
 between $2^4$ and $2^5$ & 2.08333\% \\
 between $2^5$ and $2^6$ & 0.0\% \\
 between $2^6$ and $2^7$ & 0.0\% \\
 between $2^7$ and $2^8$ & 0.10610\% \\
 between $2^8$ and $2^9$ & 0.0\% \\
 between $2^9$ and $2^{10}$ & 0.0\% \\
 between $2^{10}$ and $2^{11}$ & 0.03791\% \\
 between $2^{11}$ and $2^{12}$ & 0.01085\% \\
 between $2^{12}$ and $2^{13}$ & 0.03607\% \\
 between $2^{13}$ and $2^{14}$ & 0.00794\% \\
 between $2^{14}$ and $2^{15}$ & 0.00583\% \\
 between $2^{15}$ and $2^{16}$ & 0.14600\% \\
 between $2^{16}$ and $2^{17}$ & 0.00474\% \\
 between $2^{17}$ and $2^{18}$ & 0.00359\% \\
 between $2^{18}$ and $2^{19}$ & 0.00600\% \\
 between $2^{19}$ and $2^{20}$ & 0.00234\% \\
 between $2^{20}$ and $2^{21}$ & 0.00556\% \\
 between $2^{21}$ and $2^{22}$ & 0.00145\% \\
\end{tabular}
\end{center}

Again, the general pattern is a decrease in the rate, but the rates are not strictly decreasing.

By nature of the algorithm, successful congruence classes must have a sequence pattern that supports explicit listing of at least some of the sequence elements. To gain insight as to how far the algorithm may need to go before approaching completion, consider the sequence for the number 27. This sequence does not join a previous sequence until the 96th element when it joins the sequence for 15 at the second element (the number 46). This sequence endures 58 divisions by 2. It is reasonable to expect that if an congruence class is to be found that includes 27, the modulus will be at least $2^{58} \approx 2.8823\times 10^{17}$.

For numbers less than 10,000, the longest sequence to join a previous sequence is number 703 at element 133. Under 100,000 the longest ``modified stopping time'' is 220 elements for number 35,655. Under 1,000,000 the longest is 327 elements for number 803,871.

This crude analysis indicates a daunting numerical task ahead. The algorithm is not yet well understood. Even for the sequence for the number 27, the algorithm may travel much beyond the modulus $2^{58}$ before finding a successful congruence class containing 27. But, if an upper bound exists for the stopping time of $3n+1$ sequences, then this algorithm has the possibility of succeeding.

The list of interim results follows.

\section{Interim Results}

At the time of writing, 13,011 successful congruence classes have been identified yielding a completion percentage of approximately 99.3804\%. The LCM of the moduli is currently 2,229,025,112,064. The successful sequence patterns are given below.

Presentation of the results in the table below has been streamlined to save space. The first two columns describe the sequence pattern $bk-c$ where $b$ is the modulus and $c$ is the remainder. The third column indicates the modified stopping time of the sequence. This is either the element index at which the sequence becomes less than the anchor (in this case columns 4-6 are empty) or the index of the element that appears in a previous sequence. For the latter case, the last three columns describe the previous sequence joined and the index of the join.

Example 1: / 16 / 13 / 7 / - / - / - / would describe the sequence $16k-13$ (congruence class $\overline{13}_{16}$) which becomes less than the anchor at element 7.

Example 2: / 18 / 5 / 1 / 16 / 5 / 6 / would describe the sequence $18k-5$ (congruence class $\overline{5}_{18}$) whose 1st element is the 6th element of $16k-5$.

\begin{center}
\scriptsize
%\small
\begin{longtable}{ |r|r|r|r|r|r| }
\hline
\multicolumn{3}{|c|}{Sequence subset} & \multicolumn{3}{|c|}{Joins sequence} \\
\hline
Factor & Offset & MST & Factor & Offset & Index \\
\hline
2 & 0 & 1 & - & - & - \\
4 & 3 & 4 & - & - & - \\
6 & 1 & 1 & 4 & 1 & 3 \\
16 & 13 & 7 & - & - & - \\
18 & 5 & 1 & 16 & 5 & 6 \\
32 & 9 & 9 & - & - & - \\
\hline
\end{longtable}
\end{center}

\end{document}